\documentclass[preprint,12pt]{elsarticle}

\usepackage{graphicx}
\usepackage{amsmath}
\usepackage{amssymb}
\usepackage{amsfonts}
\usepackage{amscd}
\usepackage{latexsym}
\usepackage{graphics}
\usepackage[all]{xy}

 \newtheorem{theorem}{Theorem}[section]
  \newtheorem{prop}[theorem]{Proposition}
  \newtheorem{cor}[theorem]{Corollary}
  \newtheorem{lemma}[theorem]{Lemma}

\journal{Advances in Mathematics}

\begin{document}

\begin{frontmatter}



\title{Kato's Residue Homomorphisms and Reciprocity Laws on Arithmetic Surfaces}


\author{Dongwen Liu}

\address{Department of Mathematics, University of Connecticut, 196 Auditorium Rd, Storrs, CT 06269, USA

Email:  dongwen.liu@uconn.edu}

\begin{abstract}
We explicitly study Kato's residue homomorphisms in Milnor $K$-theory, which are closely related to Contou-Carr\`ere symbols. As applications we establish several reciprocity laws for certain locally defined maps on $K$-groups that are associated to arithmetic surfaces.
\end{abstract}


\end{frontmatter}

\section{Introduction}\label{sec1}

For a projective algebraic curve $C$ over a perfect field $k$, the famous Weil reciprocity law states that for $f, g\in k(C)^\times$, $\{f,g\}_p=1$
for almost all closed points $p\in C$, and
\[
\prod_{p\in C}N_{k(p)/k}\{f,g\}_p=1,
\]
where $k(p)$ is the residue field of $p$, $N_{k(p)/k}$ stands for the norm map, and
$\{f,g\}_p$ is the one-dimensional tame symbol defined by
\[
\{f,g\}_p=(-1)^{\nu_p(f)\nu_p(g)}\frac{f^{\nu_p(g)}}{g^{\nu_p(f)}}(p).
\]

The Weil reciprocity law, as well as the analogous reciprocity law for residues of differential forms, can be proved by reduction to
$\mathbb{P}_k^1$ (see e.g. \cite{S}), in which case everything can be calculated explicitly. On the other hand, in \cite{T} Tate gave
an intrinsic proof of the residue formula for differential forms by defining the local residues in terms of traces of certain linear operators
on infinite dimensional vector spaces. This approach and its multiplicative analog were further developed by many authors (see e.g. \cite{ADCK} \cite{APR}  \cite{OZ} \cite{PR2}) with applications to reciprocity laws on algebraic curves and surfaces. In particular, Osipov and Zhu \cite{OZ}
gave a categorical proof of the Parshin reciprocity laws for two-dimensional tame symbols on algebraic surfaces, where the notion of
non-strictly commutative Picard groupoid (cf. \cite{BBE}) was adopted.

Instead of algebraic surfaces, in this paper we will study the symbols and reciprocity laws on an arithmetic surface $X\to S$, where $S$ is the spectrum
of a characteristic zero Dedekind domain ${\mathcal O}_K$ with finite residue fields. For example we may take ${\mathcal O}_K$ to be the ring of integers of a number
field $K$. For any closed point $x$ of $X$ lying over a closed point $s$ of $S$, and an irreducible curve $y\subset X$ passing through $x$, we will obtain  a product of two-dimensional local fields $K_{x,y}$ and a local symbol $\{,\}_{x,y}$, which is a Steinberg symbol
\[
\{,\}_{x,y}: K_{x,y}^\times\times K_{x,y}^\times \to K_s^\times,
\]
where $K_s$ is the completion of $K=\textrm{Frac}(\mathcal{O}_K)$ at $s$. Then our main results concerning reciprocity laws are the following

\begin{theorem}
(i) Reciprocity around a fixed point $x$ (Theorem \ref{recipt}): Let $f, g\in K(X)^\times$. Then $\{f,g\}_{x,y}=1$ for almost all curves $y$ passing through $x$, and $\prod_{y\ni x}\{f,g\}_{x,y}=1$.

(ii) Reciprocity along a vertical curve $y$ (Theorem \ref{recivertical}): Let $f,g\in K(X)_y^\times$.  Then $\prod_{x\in y}\{f,g\}_{x,y}$ converges to 1 in $K_s^\times$, where
$K(X)_y$ is the fraction field of $\widehat{\mathcal O}_{X,y}$, the completion of the discrete valuation ring $\mathcal{O}_{X,y}$.

(iii) Reciprocity along a horizontal curve $y$ (Theorem \ref{recihori}): assume $K$ is a number field, compactify $X$ and $S$ by including archimedean places, and use an idele class character $\chi=(\chi_s)_{s\in S}$ to induce local symbols
\[
\xymatrix{
\chi_{x,y}: K_{x,y}^\times\times K_{x,y}^\times \ar[r]^{ \qquad \{,\}_{x,y}} &  K_s^\times \ar[r]^{\chi_s} &  S^1.
}
\]
 Let $f, g\in K(X)_y^\times$. Then $\chi_{x,y}(f,g)=1$ for almost all $x\in y$, and $\prod_{x\in y}\chi_{x,y}(f,g)=1$.
\end{theorem}

These results are inspired by many works, among which we would like to mention Kato's work \cite{K} on residue homomorphisms in Milnor $K$-theory, Morrow's work \cite{Mo} \cite{Mo2} on dualising sheaves and residues of differential forms on arithmetic surfaces, all the analogues of these results in the geometric case proved by Osipov \cite{O}, and more recently Osipov and Zhu's work \cite{OZ} on Parshin reciprocity laws on algebraic surfaces. Our basic settings and strategy of proofs combine certain ideas in the papers cited above, and the new ingredients will be discussed in the paragraphs below. We remark that it will be interesting to investigate the relations between residues and the multiplicative analog via certain exponential homomorphisms, see e.g. \cite{F} \cite{Ku}. Also for future works, it is natural and important to consider applications to higher adelic method, e.g. generalizing the adelic approach to intersection theory, chern classes and $L$-functions in \cite{Pa} from algebraic surfaces to arithmetic surfaces.

To construct the local symbols, given a two-dimensional local field $F$ and an embedding of a local field $k$ into $F$, we will define a Steinberg symbol
\[
\{,\}_{F/k}: F^\times\times F^\times\to k^\times.
\]
In contrast to the categorical approach as developed in \cite{OZ} for algebraic surfaces, our definition of the symbol makes use of explicit structures of two-dimensional local fields, in particular that of the so-called standard fields of mixed characteristic. In that case, we will apply Kato's residue homomorphisms in terms of Milnor $K$-theory, and write down the explicit formula
(see section \ref{3.2}) for the symbol.
The definition of the symbol extends to non-standard fields once we apply the norm map in $K$-theory and establish the functorial properties. The passage from standard fields to non-standard ones can be viewed as generalizations of the classical methods of ``reduction to $\mathbb{P}^1$", which we mentioned at the beginning.

An important feature of Kato's residue symbol for standard fields is that, modulo a power of the maximal ideal it coincides with the famous Contou-Carr\`ere symbol (see section \ref{3.3}). This fact was proved by P\'al in \cite{P}. Thus we may derive the rigidity of the symbol from the fact that Contou-Carr\`ere symbols are invariant under reparametrizations. On the other hand, this connection suggests that it is also likely to give an intrinsic treatment of Tate's style (cf. \cite{APR}), which  might involve functional analysis on two-dimensional local fields \cite{Ca}. We hope to address this question somewhere else.

The paper is organized as follows. In section \ref{sec2} as preliminaries we briefly review Kato's residue homomorphisms in Milnor's $K$-theory. In section \ref{sec3} we apply Kato's theory to define the symbols for two-dimensional local fields, establish the rigidities and functorialities, and discuss the relations with Contou-Carr\`ere symbols. In section \ref{sec4} we apply previous results to prove a reciprocity for certain two-dimensional normal local rings, and then geometrize to prove the reciprocity around a closed point on an arithmetic surface. In section \ref{sec5} we prove the reciprocities along irreducible vertical and horizontal curves, using certain density arguments. The treatment in the last two sections follows \cite{Mo2} closely. Analogous counterparts of many propositions and proofs can be found in \cite{Mo2}.


\section{Milnor's $K$-theory and Kato's residue homomorphisms}\label{sec2}

In this section we briefly review Kato's definition of residue homomorphisms for two-dimensional local fields of mixed characteristic.

Throughout this paper, we use $K_n$ to denote Milnor, rather than Quillen $K$-theory, with the only exception in section 3.3 where we will consider Quillen $K$-groups for semi-local rings. Although for our purpose we are only concerned with $K_2$, let us recall the general construction from \cite{K} for completeness.
For a field $k$, let $K_n(k)$ ($n\geq 0$) be the $n$th Milnor's $K$-group of $k$ (cf. \cite{M} \cite{M2}), i.e.
$K_n(k)=(k^\times)^{\otimes n}/J,$
where $J=\langle a_1\otimes a_2\otimes\cdots\otimes a_n, a_i+a_j=1\textrm{ for some }i\neq j\rangle$.
For a discrete valuation field $k$, let $\nu_k$ be the normalized additive discrete valuation of $k$. Let
\begin{eqnarray*}
&&{\mathcal O}_k=\{x\in k: \nu_k(x)\geq 0\},\quad \mathfrak{m}_k=\{x\in k: \nu_k(x)\geq 1\},\\
&&U^i_k=\{x\in k: \nu_k(x-1)\geq i\},~ i\geq 1.
\end{eqnarray*}
The residue class field of $k$ is denoted by $\bar{k}$, and for $x\in {\mathcal O}_k$ let $\bar{x}$ be the residue class of $x$ in $\bar{k}$. For $i\geq 1$ let $\bigoplus_{n\geq 0} U^i K_n(k)$ be the graded ideal of $\bigoplus_{n\geq 0}K_n(k)$ generated by elements of $U_k^i\subset k^\times= K_1(k)$. Let
\[
\widehat{K}_n(k)=\varprojlim_i K_n(k)/U^i K_n(k).
\]
If $\hat{k}$ is the completion of $k$, then there are canonical isomorphisms $K_n(k)/U^iK_n(k)\cong K_n(\hat{k})/U^i K_n(\hat{k})$, $i\geq 1$ and $\widehat{K}_n(k)\cong \widehat{K}_n(\hat{k})$.

Let us recall the boundary homomorphism $\partial$ and the norm homomorphism of Milnor's $K$-group \cite{BT}. For a discrete valuation field $k$, there exists the \textit{higher tame symbol}, which is the unique
homomorphism
\[
\partial: K_{n+1}(k)\to K_n(\bar{k})\quad (n\geq 0)
\]
such that
\[
\partial\{x_1,\cdots, x_n, y\}=\nu_k(y)\cdot\{\bar{x}_1,\cdots, \bar{x}_n\}
\]
for any $x_1,\cdots, x_n\in {\mathcal O}_k^\times$ and $y\in k^\times$. For a given parameter $\pi$, there is also a \textit{specilization map}
\[
\lambda_\pi: K_n(k)\to K_n(\bar{k})
\]
defined by
\[
\lambda_\pi\{u_1\pi^{i_1},\ldots, u_n\pi^{i_n}\}=\{\bar{u}_1,\ldots, \bar{u}_n\}
\]
where $u_1,\ldots, u_n\in {\mathcal O}^\times_k$.

For any finite field extension $E/F$, there exists a canonical norm homomorphism, or \textit{transfer map}, $N_{E/F}: K_n(E)\to K_n(F)$, which satisfies the projection formula, functoriality, and reciprocity (cf. \cite{BT}).
Let $k'$ be a finite extension of a discrete valuation field $k$. If the integral closure of ${\mathcal O}_k$ in $k'$ is a finitely generated ${\mathcal O}_k$-module (e.g. if $k$ is complete), then the following diagram is commutative
\begin{equation}\label{ext}
\xymatrix{
K_{n+1}(k') \ar[d]_{N_{k'/k}} \ar[r]^{(\partial_v)_v} & \bigoplus\limits_v K_n(\kappa(v)) \ar[d]^{\sum_v N_{\kappa(v)/\bar{k}}} \\
K_{n+1}(k) \ar[r]^{\partial} & K_n(\bar{k})
}
\end{equation}
where $v$ runs over all normalized discrete valuations of $k'$ such that $\{x\in k: v(x)\geq 0\}={\mathcal O}_k$, and $\kappa(v)$ is the residue field of
$v$.

Now let us give Kato's definition of the residue homomorphism. Let $k$ be a complete discrete valuation field, and let $M$ be the fraction field
of ${\mathcal O}_k[[T]]$. Let $k\{\{T\}\}$ be the completion of $M$ with respect to the discrete valuation of $M$ defined by the height 1 prime ideal
$\mathfrak{m}_k{\mathcal O}_k[[T]]$ of ${\mathcal O}_k[[T]]$. In concrete terms, $k\{\{T\}\}$ is the field of all formal series $\sum\limits_{i\in\mathbb{Z}}a_i T^i$
over $k$ such that $\nu_k(a_i)$ is bounded below and $\lim\limits_{i\to -\infty}a_i=0$. The valuation $\nu_{k\{\{T\}\}}$ is given by $\inf\limits_{i\in\mathbb{Z}}\nu_k(a_i)$ and the residue field of $k\{\{T\}\}$ is $\bar{k}((T))$.

Let $\mathfrak{S}=\{\mathfrak{p}: \mathfrak{p}\textrm{ is a height 1 prime ideal of }{\mathcal O}_k[[T]], \mathfrak{p}\neq\mathfrak{m}_k{\mathcal O}_k[[T]]\}$. By Weierstrass's preparation theorem, each element of $\mathfrak{S}$ is generated by an irreducible, distinguished polynomial (i.e. of the form $T^l+a_1T^{l-1}+\cdots+a_l$ with $a_i\in\mathfrak{m}_k$). Let $res_M: K_{n+1}(M)\to K_n(k)$ be the composition map
\begin{equation}\label{kato}
K_{n+1}(M) \stackrel{(\partial_{\mathfrak{p}})_{\mathfrak{p}}}{\longrightarrow}  \bigoplus_{\mathfrak{p}\in\mathfrak{S}}K_n(\kappa(\mathfrak{p})) \stackrel{\sum_{\frak{p}}N_{\kappa(\mathfrak{p})/k}}{\longrightarrow} K_n(k).
\end{equation}
In \cite{K} K. Kato proved the following

\begin{theorem}\label{Kato}
The homomorphism $res_M: K_{n+1}(M)\to K_n(k)$ satisfies
\[
res_M(U^iK_{n+1}(M))\subset U^iK_n(k) \textrm{ for any }i\geq 1,
\]
and hence induces a homomorphism
\[
res_{k\{\{T\}\}}:~\widehat{K}_{n+1}(k\{\{T\}\})\cong \widehat{K}_{n+1}(M)\to\widehat{K}_n(k).
\]
\end{theorem}

For any field extension $E/F$ let $j_{E/F}$ stand for the natural map $K_n(F)\to K_n(E)$. We have the functoriality $j_{E/E}=\rm{Id}$ and $j_{E/F}\circ j_{F/K}=j_{E/K}$. If $E/F$ has degree $d$, then the composition
\begin{equation*}
\xymatrix{
K_n(F)\ar[r]^{j_{E/F}} & K_n(E)\ar[r]^{N_{E/F}} & K_n(F)
}
\end{equation*}
is multiplication by $d$. Note that in general $j_{E/F}$ is not injective, but its kernel is torsion provided that $E$ is algebraic over $F$.

The following local-global norm formula for $K$-groups is standard, see \cite{BT} and \cite{GS} Lemma 7.3.6.

\begin{prop}\label{lgnorm}
Let $E/F$ be a separable field extension. Then for any discrete valuation $y$ of $F$, and any $n\geq 0$, we have
\[
j_{F_{y}/F}\circ N_{E/F}=\prod_{Y|y}N_{E_{Y}/F_{y}}\circ j_{E_{Y}/E}: K_n(E)\to K_n(F_y),
\]
where $Y$ runs over all discrete valuations of $E$ which extend $y$.
\end{prop}

Finally let us recall a very useful theorem of Kato (\cite{K} Proposition 2), which states that if $k'/k$ is a finite extension of complete discrete valuation fields, then $N_{k'/k}U^{ie}K_n(k')\subset U^i K_n(k)$ for all $i\geq 1$, where $e=e(k'/k)$ is the ramification index. This result was recently reproved in an explicit way by M. Morrow in \cite{Mo3}.

\section{Symbols on two-dimensional local fields}\label{sec3}

Most of the materials in this section can be viewed as a survey of Kato's residue symbol \cite{K} in case of two-dimensional local fields and the Contou-Carr\`ere symbol \cite{CC}. We will give some explicit descriptions and then discuss the relation between these two symbols. Moreover, certain expected functoriality results will be proved.

\subsection{Classical tame symbols}\label{3.1}

From now on, local fields are always assumed to have finite residue fields. A two-dimensional local field is a complete discrete valuation field $F$ whose residue field $\overline{F}$ is a local field. We are only interested in the case that $F$ is of characteristic zero and $\overline{F}$ is non-archimedean. Then we say that $F$ has equal or mixed characteristic according to $\overline{F}$
has characteristic zero or not. Assume that an embedding of a local field $k$ into $F$ has been fixed. We will define a Steinberg symbol
\[
\{,\}_{F/k}: F^\times\times F^\times\to k^\times.
\]
And we shall denote by $\partial_{F/k}$ the induced homomorphism
$\partial_{F/k}: K_2(F)\to K_1(k)$.

In the equal characteristic case, let $k_F$ be the algebraic closure of $k$ inside $F$, which is called the \textit{constant field} of $F$. Then $k_F$ is a finite extension of $k$ and there is
a $k_F$-isomorphism $F\cong k_F((t))$ by choosing a uniformiser $t\in F$. We define $\{,\}_{F/k}$ to be the composition
\begin{equation}\label{equtame}
\xymatrix{
F^\times\times F^\times \ar[r]^{\{,\}_F} & k_F^\times \ar[r]^{N_{k_F/k}} & k^\times,
}
\end{equation}
where $\{,\}_F$ is the classical tame symbol
\begin{equation}\label{classtame}
\{f,g\}_F=(-1)^{\nu_F(f)\nu_F(g)}\frac{f^{\nu_F(g)}}{g^{\nu_F(f)}}\mod \mathfrak{m}_F.
\end{equation}
Note that in this case $\overline{F}=k_F$, and we shall denote by $\partial_F$ the boundary map $\partial: K_2(F)\to K_1(k_F)$. The following functoriality  can be deduced as a corollary of
(\ref{ext}).

\begin{prop}\label{funequ}
Let $F'$ be a finite extension of $F$. Then the following diagram commutes:
\begin{equation*}
\xymatrix{
K_2(F') \ar[d]_{N_{F'/F}} \ar[r]^{\partial_{F'}} & K_1(k_{F'})\ar[d]^{N_{k_{F'}/k_F}}\\
K_2(F) \ar[r]^{\partial_F} & K_1(k_F).
}
\end{equation*}
\end{prop}

\subsection{Kato's residue homomorphisms for standard fields}\label{3.2}

From now on assume that $\overline{F}$ has characteristic $p$. Let $k_F$ be the algebraic closure of $\mathbb{Q}_p$ inside $F$, which is called
the \textit{coefficient field} of $F$. Then
$k\subset k_F$ are finite extensions of $\mathbb{Q}_p$. Similarly as the equal characteristic case, we will define a symbol
\begin{equation}\label{mixtame}
\{,\}_F: F^\times\times F^\times\to k_F^\times
\end{equation}
and compose with the norm map $N_{k_F/k}$ to obtain $\{,\}_{F/k}$.

By structure theory of two-dimensional local fields of mixed characteristic, there is a two-dimensional local field $L$ inside $F$ such that $F/L$ is finite, $\overline{L}=\overline{F}$, $k_L=k_F$ and $L$ is $k_L$-isomorphic to $k_L\{\{T\}\}$, see \cite{C} and \cite{Mo} Lemma 2.14.
With the last property, $L$ is called a \textit{standard} two-dimensional local field of mixed characteristic.  We will study only standard field $L$ in this section and the next, then return to $F$ in section 3.4.

Let $M$ be the field of fractions of $A={\mathcal O}_{k_L}[[T]]$. Due to Weierstrass's preparation theorem (see e.g \cite{W} Theorem 7.3), any $f\in M^\times$ can be uniquely factored as
\begin{equation}
 f=f_0\cdot\dfrac{a_f}{b_f}\cdot u_f,
 \end{equation}
where $f_0\in k_L^\times$, $a_f$, $b_f\in{\mathcal O}_{k_L}[T]$ are relatively prime distinguished polynomials, and $u_f\in 1+T{\mathcal O}_{k_L}[[T]]$.
 Define the \textit{winding number} of $f$ to be
 \begin{equation}\label{wind}
 w(f)=\deg(a_f)-\deg(b_f).
 \end{equation}
 Fix an algebraic closure $k_L^{al}$ of $k_L$. For any $p(T)\in k_L[T]$ denote by $V(p)$ the set of roots of $p$ in $k_L^{al}$ counting multiplicities.
We now explicitly describe Kato's residue symbol, by defining a map $\{,\}_M: M^\times\times M^\times \to k_L^\times$:
\begin{equation}\label{tame}
\{f, g\}_M=(-1)^{w(f)w(g)}\frac{{g_0}^{w(f)}}{f_0^{w(g)}}\frac{\prod\limits_{\alpha_f\in V(a_f)}u_g(\alpha_f)}{\prod\limits_{\alpha_g\in V(a_g)}u_f(\alpha_g)}
\frac{\prod\limits_{\beta_g\in V(b_g)}u_f(\beta_g)}{\prod\limits_{\beta_f\in V(b_f)}u_g(\beta_f)}.
\end{equation}
Since roots of a distinguished polynomial lie in $\mathfrak{m}_{k_L^{al}}$, the right-hand-side of (\ref{tame}) makes sense and belongs to $k_L^\times$,
which can be easily seen.

We will see shortly (Proposition \ref{=res} below) from P\'al's work \cite{P} that $\{,\}_M$ is induced from the negative of Kato's residual symbol $res_M$.
Before that let us first establish some basic properties of $\{,\}_M$ and extend the symbol to the standard field $L$.

\begin{lemma}
$\{,\}_M$ is a Steinberg symbol, i.e. is bi-multiplicative and satisfies the Steinberg property $\{f,1-f\}_M=1$ for $f\in M^\times$, $f\neq 1$.
\end{lemma}

\noindent\textit{Proof.}
$\{,\}_M$ is obviously bi-multiplicative. Let us prove the Steinberg property, which is in fact an elementary exercise. Write $g=1-f$ and suppose that $f, g$ decompose as above. Then it is easy to see that $b_f=b_g$. Since $f+g=1$, for any $\alpha_f\in V(a_f)$ one has $g(\alpha_f)=1$, which implies that
\[
u_g(\alpha_f)=g_0^{-1}\frac{b_g(\alpha_f)}{a_g(\alpha_f)}.
\]
We also have a similar equation with $f$ and $g$ exchanged. For any $\beta\in V(b_f)=V(b_g)$, again from $f+g=1$ we may derive that
\[
f_0a_f(\beta)u_f(\beta)=-g_0 a_g(\beta)u_g(\beta), \quad \frac{u_f(\beta)}{u_g(\beta)}=-\frac{g_0}{f_0}\frac{a_g(\beta)}{a_f(\beta)}.
\]
For two monic polynomials $p(T), q(T)\in k_L[T]$, by factoring into linear factors one may deduce the formula
\[
\prod_{\alpha\in V(p)}q(\alpha)=(-1)^{\deg p\cdot\deg q}\prod_{\beta\in V(q)}p(\beta).
\]
Applying this formula for $(p, q)=(a_f, a_g)$, $(a_f, b_f)$ and $(a_g, b_g)$, from all the above one can check that $\{f,g\}_M=1$.
\hfill$\Box$

\

Therefore $\{,\}_M$ descends to a homomorphism $\partial_M: K_2(M)\to K_1(k_L)$. We have the following

\begin{prop}\label{cont}
The homomorphism $\partial_M: K_2(M)\to K_1(k_L)$ satisfies
\[
\partial_M(U^iK_2(M))\subset U^iK_1(k_L)=U^i_{k_L} \textrm{ for any }i\geq 1,
\]
hence induces a homomorphism
\[
\widehat{\partial}_L: \widehat{K}_2(L) \cong \widehat{K}_2(M)\to \widehat{K}_1(k_L)=K_1(k_L),
\]
noting that $k_L$ is complete.
\end{prop}

\noindent\textit{Proof.}
Let $\mathfrak{B}$ be the set of all normalized additive discrete valuations $v$ of $k_L(X)$ such that
  $v(k_L^\times)=0$, i.e. corresponding to the set of closed points of $\mathbb{P}^1_{k_L}$. Let $\mathfrak{S}$ be the subset
  of $\mathfrak{B}$ corresponding to irreducible distinguished polynomials in ${\mathcal O}_{k_L}[T]$. Let $A_+={\mathcal O}_{k_L}[[T]]\otimes_{{\mathcal O}_{k_L}}k_L$,
$A_-=\{f\in k_L(X): v(f)\geq 0 \textrm{ if }v\in \mathfrak{B}-\mathfrak{S}-\{\infty\}\}$. By \cite{K} Proposition 1, $U^iK_n(M)=U^iK_n(A_+)+U^iK_n(A_-)$ (see [\textit{loc. cit.}] for the precise definitions of these notations). From the definition of $\{,\}_M$ it is clear that
$\partial_M(U^iK_2(A_+))=1$ for any $i\geq 1$.

Let us prove that $\partial_M(U^iK_2(A_-))\subset U^iK_1(k_L)$ for any $i\geq 1$, i.e. prove that
$\{f,g\}_M\in U^i_{k_L}$ for any $f\in A_-^\times$, $g\in U^iK_1(A_-)$. We may write $g$ as $g_0\dfrac{a_g}{b_g}$ such that $g_0\in U^i_{k_L}$, $a_g$ and $b_g$ are distinguished polynomials such that $a_g\equiv b_g\mod \mathfrak{m}_{k_L}^i{\mathcal O}_{k_L}[T]$. In particular $a_g$ and $b_g$ have the same degree, i.e. $w(g)=0$. Then we have
\[
\{f, g\}_M={g_0}^{w(f)}\frac{\prod_{\beta \in V(b_g)}u_f(\beta)}{\prod_{\alpha\in V(a_g)}u_f(\alpha)}.
\]
From $a_g\equiv b_g\mod \mathfrak{m}_{k_L}^i{\mathcal O}_{k_L}[T]$, it is easy to show that $\prod_{\beta \in V(b_g)}u_f(\beta)\equiv \prod_{\alpha\in V(a_g)}u_f(\alpha)\mod \mathfrak{m}_{k_L}^i$. Indeed, let $n=\deg(a_g)=\deg(b_g)$ and consider the series $\prod^n_{j=1}u_f(x_j)$ in the variables $x_1,\ldots,x_n$. Then one may write the series as a formal sum of symmetric polynomials in $x_1,\ldots,x_n$. If we substitute $\alpha\in V(a_g)$ or
$\beta\in V(b_g)$ for these $x_i$'s, then the values of symmetric polynomials can expressed in terms of the coefficients in $a_g$ or $b_g$. The convergence of the formal sum is obvious. Hence from above arguments it follows that $\{f,g\}_M\in U^i_{k_L}$.
\hfill$\Box$

\

Now define $\{,\}_L$ to be the composition
\begin{equation}\label{std}
\xymatrix{
L^\times\times L^\times \ar[r] & K_2(L) \ar[r] & \widehat{K}_2(L) \ar[r]^{\widehat{\partial}_L} & K_1(k_L),
}
\end{equation}
and $\{,\}_{L/k}=N_{k_L/k}\circ\{,\}_L$. We also denote by $\partial_L$ the map $K_2(L)\to K_1(k_L)$ in (\ref{std}). Before we proceed further let us give the explicit formula of $\{,\}_L$ (see \cite{A} section 4 for a similar treatment). Any $f\in L^\times$ can be expressed in a unique way (see e.g. \cite{A} section 4.1) as
\begin{equation}\label{witt}
f=f_0T^{w(f)}\prod^\infty_{i=1}(1-f_{-i}T^{-i})(1-f_iT^i),
\end{equation}
with $w(f)\in\mathbb{Z}$, $f_0\in k_L^\times$, $f_i\in {\mathcal O}_{k_L}$, $f_{-i}\in\mathfrak{m}_{k_L}$ for $i>0$, and $\lim\limits_{i\to\infty}f_{-i}=0$. In (\ref{witt}) we call $w(f)$ the \textit{winding number} of $f$, extending the definition (\ref{wind}), and we call the sequence $\{f_i\}_{i\in\mathbb{Z}}$ the \textit{Witt parameters} of $f$.

\begin{lemma}
If $f, g\in L^\times$ have winding numbers $w(f), w(g)$ and Witt parameters $\{f_i\}, \{g_j\}$ respectively, then
\begin{equation}\label{tame2}
\{f,g\}_L=(-1)^{w(f)w(g)}\frac{g_0^{w(f)}}{f_0^{w(g)}}\prod_{i,j\geq 1}\frac{(1-f_{-i}^{j/(i,j)}g_j^{i/(i,j)})^{(i,j)}}{(1-f_i^{j/(i,j)}g_{-j}^{i/(i,j)})^{(i,j)}}
\end{equation}where $(i,j)$ is the greatest common divisor of $i$ and $j$.
\end{lemma}

\noindent\textit{Proof.}
Consider the special case $f=1-a T^{-i}$, $g=1-bT^j$ with $a\in \frak{m}_{k_L},$ $b\in\mathcal{O}_{k_L}$. Let $\alpha\in k_L^{al}$ be an $i$th root of $a$, and $\xi\in k_L^{al}$ be a primitive $i$th root of unity. Then
$\xi^j$ is a primitive $i/(i,j)$-th root of unity, and thus one has
\[
\prod^i_{l=1}(1-b\alpha^j\xi^{jl})=(1-b^{i/(i,j)}\alpha^{ji/(i,j)})^{(i,j)}=(1-b^{i/(i,j)}a^{j/(i,j)})^{(i,j)}.
\]
Thus (\ref{tame2}) holds in this case. In a similar fashion, it is also straightforward to check other cases, e.g. $f=1-aT^i$, $g=1-bT^j$ with $a, b\in{\cal O}_{k_L}$,
or $f=1-aT^{-i}$, $g=1-bT^{-j}$ with $a, b\in\frak{m}_{k_L}$ (see also \cite{Pa} Lemmas 3.4 and 3.5).

Now let us consider the general case. For $N\geq 0$ let
\begin{eqnarray*}
&&f_N^+=\prod^N_{i=1}(1-f_{i}T^{i}),\quad f_N^-=\prod^N_{i=1}(1-f_{-i}T^{-i}), \quad f^{(N)}=f_0T^{w(f)}f_N^+f_N^-,\\
 &&\tilde{f}_N^+=\prod_{i>N} (1-f_iT^i),\quad \tilde{f}_N^-=\prod_{i>N}(1-f_{-i}T^{-i}).
\end{eqnarray*}
Then $f=f^{(N)}\tilde{f}_N^+\tilde{f}_N^-$. Define $g^{(N)}$, $g_N^\pm$, $\tilde{g}_N^\pm$ similarly.
Then we have
\[
\{f^{(N)}, g^{(N)}\}_L=(-1)^{w(f)w(g)}\frac{g_0^{w(f)}}{f_0^{w(g)}}\prod_{i,j=1}^N\frac{(1-f_{-i}^{j/(i,j)}g_j^{i/(i,j)})^{(i,j)}}{(1-f_i^{j/(i,j)}g_{-j}^{i/(i,j)})^{(i,j)}}
\]
Note that the right-hand-side of (\ref{tame2}) is a convergent product in $k_L^\times$, we only need to show that $\{f^{(N)}, g^{(N)}\}_L\to \{f,g\}_L$ as $N\to\infty$.
 From $\lim\limits_{i\to\infty}f_{-i}=0$ it follows that $\tilde{f}_N^-\to 1$ and $f_N^-\to \tilde{f}_0^-$ as $N\to\infty$. Since $\{,\}_L$ is continuous by Proposition \ref{cont}, it suffices to prove that $\{\tilde{f}_N^+,g\}_L\to 1$. From (\ref{tame}) we obtain $\{\tilde{f}_N^+, g_0T^{w(g)}\tilde{g}^+_0\}_L=1$ and
 \[
 \{\tilde{f}_N^+,g\}_L= \{\tilde{f}_N^+,\tilde{g}_0^-\}_L =\lim_{j\to\infty} \{\tilde{f}_N^+, g_j^-\}_L=\prod_{i>N, j>0}(1-f_i^{j/(i,j)}g_{-j}^{i/(i,j)})^{-1},
 \]
 which converges to $1$ as $N\to\infty$. In the above we used $g_j^-\to \tilde{g}_0^-$ and again the continuity of $\{,\}_L$, while $\{\tilde{f}_N^+, g_j^-\}_L$ was calculated using (\ref{tame}), similar to the special case we considered at the beginning of the proof. This justifies the lemma.
\hfill$\Box$

\begin{prop}\label{=res}
We have $\widehat{\partial}_L=-res_L: \widehat{K}_2(L)\to K_1(k_L)$, where $\widehat{\partial}_L$ and $res_L$ are defined by Proposition \ref{cont} and Theorem \ref{Kato} respectively.
\end{prop}

\noindent \textit{Proof.} The assertion follows from \cite{P} Lemmas 3.4 and 3.5 immediately, noting that both symbols are continuous.
\hfill$\Box$

\

We remark that given Proposition \ref{=res}, a priori we will have the continuity of $\{,\}_L$, and thus Proposition \ref{cont} is an immediate consequence.

\subsection{Rigidity of Kato's residue symbol: reduction to Contou-Carr\`{e}re symbols}\label{3.3}

In \cite{K} Kato proved the rigidity of the residue homomorphism, which implies the independence of the choice of the local parameter $T$. In conjunction with Proposition \ref{=res}, this guarantees the well-definedness of the symbol $\{,\}_L$ for standard fields. Since Kato's proof involves a certain amount of labors, for the case of standard fields we offer an alternative proof which meanwhile serves another purpose of illustrating the connections with the Contou-Carr\`{e}re symbols. Similar results are given by P\'al in \cite{P} section 4.

For completeness let us recall from \cite{CC} the basic theory of Contou-Carr\`{e}re symbols.
For an artinian local ring $R$ with maximal ideal $\mathfrak{m}$, any $f\in R((T))^\times$ can be written in exactly one way as
\[
f=f_0T^{w(f)}\prod^\infty_{i=1}(1-f_{-i}T^{-i})(1-f_iT^i)
\]
with $w(f)\in\mathbb{Z}$, $f_0\in R^\times$, $f_i\in R$, $f_{-i}\in\mathfrak{m}$ for $i>0$, and $f_{-i}=0$ for $i\gg 0$. We call $w(f)$ the winding number of $f$
and $\{f_i\}_{i\in\mathbb{Z}}$ the Witt parameters of $f$. The Contou-Carr\`{e}re symbol $\langle,\rangle_{R((T))}: R((T))^\times\times R((T))^\times\to R^\times$ is defined by
\begin{equation}\label{CC}
\langle f, g\rangle_{R((T))^\times}=(-1)^{w(f)w(g)}\frac{f_0^{w(g)}}{g_0^{w(f)}}\prod_{i,j\geq 1}\frac{(1-f_{i}^{j/(i,j)}g_{-j}^{i/(i,j)})^{(i,j)}}{(1-f_{-i}^{j/(i,j)}g_{j}^{i/(i,j)})^{(i,j)}}.
\end{equation}
The definition makes sense because only finitely many terms in the infinite product differ from 1. The Contou-Carr\`{e}re symbol is a Steinberg symbol, which generalizes the classical tame symbol, and also contains the residue as a special case (see the introduction of \cite{APR} for more precise statements, and \cite{PR} for a proof of the Steinberg property).

The relations between Kato's residue symbol for standard fields and Contou-Carr\`{e}re symbols are obtained in the following lemma, via reduction modulo a power of $\mathfrak{m}_{k_L}$. Let ${\mathcal O}_{k_L}\{\{T\}\}\stackrel{\rm{def}}{=}{\mathcal O}_L$ be the ring of integers of $L$, and
let $R_n={\mathcal O}_{k_L}/\mathfrak{m}_{k_L}^n$ which is an artinian local ring for $n\geq 1$. Then there is a compatible family of natural maps
${\mathcal O}_{k_L}\{\{T\}\}^\times\to R_n((T))^\times$ for $n\geq 1$ such that
\[
{\mathcal O}_{k_L}\{\{T\}\}^\times=\varprojlim_n R_n((T))^\times.
\]

\begin{lemma}\label{red}
For any $n\geq 1$, the following diagram commutes:
\begin{equation*}
\xymatrix{
{\mathcal O}_{k_L}\{\{T\}\}^\times\times{\mathcal O}_{k_L}\{\{T\}\}^\times  \ar[r] \ar[d]_{\{,\}_L} & R_n((T))^\times\times  R_n((T))^\times \ar[d]^{\langle,\rangle_{ R_n((T))^\times}^{-1}}\\
{\mathcal O}_{k_L}^\times \ar[r] & R_n^\times
}
\end{equation*}
where the left vertical arrow is the restriction of the symbol $\{,\}_L$, and the right vertical arrow is the inverse of the Contou-Carr\`ere symbol $\langle,\rangle_{R_n((T))^\times}.$
\end{lemma}

\noindent\textit{Proof.}
Direct consequences of (\ref{tame2}) and (\ref{CC}).
\hfill$\Box$

\begin{cor}
The definition of $\{,\}_L$ does not depend on the choice of the $k_L$-isomorphism $L\cong k_L\{\{T\}\}$ $($i.e. the choice
of the local parameter $T)$.
\end{cor}

\noindent\textit{Proof.}
It is known that the Contou-Carr\`ere symbol is invariant under reparametrization of $R((T))$ in the sense that, if $t\in R((T))$
is an element with winding number equal to $1$, then $\langle f, g\rangle_{R((T))^\times}=\langle f\circ t, g\circ t\rangle_{R((T))^\times}$. See \cite{PR}
for a proof of this fact, which is an application of the ``adjunction formula" [\textit{loc. cit.}] (2.3).

From this fact together with Lemma \ref{red}, by taking limit we deduce that $\{,\}_L$ is well-defined when restricted on ${\mathcal O}_L^\times \times {\mathcal O}_L^\times$. To complete the proof, we note that $L^\times={\mathcal O}_L^\times\otimes_{{\mathcal O}_{k_L}^\times}k_L^\times$, while
$\{f, a\}_L=a^{w(f)}$
for any $f\in{\mathcal O}_L^\times$, $a\in k_L^\times$. Then it remains to show that the winding number $w(f)$ is independent of the choice of $T$. In fact this is true for any $f\in L^\times$ and the assumption that $f\in\mathcal{O}_L^\times$ makes no loss of generality. From either (\ref{wind}) or (\ref{witt}) we can see that
$w(f)=\nu(\bar{f}),$ where $\bar{f}$ is the image of $f$ under the map $\mathcal{O}_L^\times\to \overline{L}^\times\cong\bar{k}_L((T))$, and $\nu$ is the normalized valuation of $\bar{k}_L((T))$. A different choice of $T$ gives another local parameter of $\bar{k}_L((T))$, hence does not change the valuation $\nu$.
The proof is finished.
\hfill$\Box$

\

Let us give the relation between Kato's residue symbol for standard fields and the tame symbol for the residue fields. We consider Quillen $K$-groups
 of a commutative semi-local ring $R$. It is known that (cf. \cite{DS})
$K_1(R)=R^\times$ and $K_2(R)$ is generated by Steinberg symbols $\{r,s\}$ with $r, s\in R^\times$, subject to certain relations. Using this it is not hard to show that by restriction one may obtain a symbol $\partial_L: K_2({\mathcal O}_L)\to K_1({\mathcal O}_{k_L})$. On the other hand, we have a $\bar{k}_L$-isomorphism $\overline{L}\cong\bar{k}_L((T))$, and thus the classical tame symbol $\partial_{\overline{L}}: K_2(\overline{L})\to K_1(\bar{k}_L)$.

\begin{cor}\label{residue}
The following diagram commutes:
\begin{equation*}
\xymatrix{
K_2({\mathcal O}_L)\ar[rr]^{\partial_L} \ar[d] && K_1({\mathcal O}_{k_L}) \ar[d] \\
K_2(\overline{L}) \ar[rr]^{\partial^{-1}_{\overline{L}}} && K_1(\bar{k}_L).
}
\end{equation*}
\end{cor}

\noindent\textit{Proof.}
This follows from Lemma \ref{red} with $n=1$. Notice that the Contou-Carr\`ere symbol
reduces to the classical tame symbol when the base artinian local ring $R$ is a field.
\hfill$\Box$

\begin{cor}
Let $k\subset L$ be a local field. Then the following diagram commutes:
\begin{equation*}
\xymatrix{
K_2({\mathcal O}_L) \ar[rr]^{\partial_{L/k}} \ar[d] && K_1({\mathcal O}_k) \ar[d]\\
K_2(\overline{L}) \ar[rr]^{e(L/k)N_{\bar{k}_L/\bar{k}}\partial^{-1}_{\overline{L}}} && K_1(\bar{k}),
}
\end{equation*}
where $e(L/k)=e(k_L/k)$ is the ramification index.
\end{cor}

\noindent\textit{Proof.}
Combine Corollary \ref{residue} with the following commutative diagram
\begin{equation*}
\xymatrix{
{\mathcal O}_{k_L}^\times \ar[rr]^{N_{k_L/k}} \ar[d] && {\mathcal O}_k^\times \ar[d]\\
\bar{k}_L^\times \ar[rr]^{e(k_L/k)N_{\bar{k}_L/\bar{k}}} && \bar{k}^\times.
}
\end{equation*}
\hfill$\Box$

Let us establish the functoriality of Kato's residue symbol for standard fields.

\begin{prop}\label{funstd}
Suppose that $L'$ is a finite extension of $L$, and that $L'$ is also standard. Then the following diagram commutes:
\begin{equation*}
\xymatrix{
K_2(L') \ar[r]^{\partial_{L'}} \ar[d]_{N_{L'/L}} & K_1(k_{L'}) \ar[d]^{N_{k_{L'}/k_L}}\\
K_2(L) \ar[r]^{\partial_L} & K_1(k_L).
}
\end{equation*}
\end{prop}

\noindent\textit{Proof.} Similar to \cite{Mo} Proposition 2.20, using the intermediate extension $Lk_{L'}$ and functoriality of norm homomorphisms, we are reduced to two cases: (i) we have compatible isomorphisms $L\cong k_L\{\{T\}\}$, $L'\cong k_{L'}\{\{T\}\}$, (ii) we have $k_{L'}=k_L$. For case (ii), by
 the usual ``principle of prolongation of algebraic identities" trick \cite{S} II.13, we may further reduce to the case $L\cong k_L\{\{T\}\}$, $L'\cong k_L\{\{t\}\}$ with $T=t^e$.

 Let $A'={\mathcal O}_{k_L'}[[T]]$ or ${\mathcal O}_{k_L}[[t]]$ for these two cases respectively, and let $M'$ be the fraction field of $A'$. Recall that previously we defined $M$ to be the fraction field of $A={\mathcal O}_{k_L}[[T]]$. Let $y$ be the height 1 prime ideal $\mathfrak{m}_{k_L}{\mathcal O}_{k_L}[[T]]$ of $A$. In both cases, $L'=M'_Y$ where $Y$ is the unique height 1 prime ideal of $A'$ lying over $y$. Applying Proposition \ref{lgnorm} for the Dedekind domains $A_y$ and $A'_y=(A\backslash y)^{-1}A'$, we obtain the following commutative diagram
\begin{equation*}
\xymatrix{
K_2(M') \ar[r]^{j_{L'/M'}} \ar[d]_{N_{M'/M}} & K_2(L')\ar[d]^{N_{L'/L}} \\
K_2(M) \ar[r]^{j_{L/M}} & K_2(L).
}
\end{equation*}
By \cite{K} Proposition 2, the norm homomorphism in Milnor $K$-theory is continuous for complete discrete valuation fields, hence passing to the limit yields
\begin{equation*}
\xymatrix{
\widehat{K}_2(M')\ar[r]^{\approx} \ar[d]_{N_{M'/M}} & \widehat{K}_2(L')\ar[d]^{N_{L'/L}}\\
\widehat{K}_2(M) \ar[r]^{\approx} & \widehat{K}_2(L).
}
\end{equation*}
To finish the proof, it remains to show that the following diagram commutes
\begin{equation*}
\xymatrix{
K_2(M') \ar[r]^{\partial_{M'}} \ar[d]_{N_{M'/M}} & K_1(k_{L'})\ar[d]^{N_{k_{L'}/k_L}} \\
K_2(M) \ar[r]^{\partial_M} & K_1(k_L).
}
\end{equation*}
This follows from (\ref{ext}), (\ref{kato}) and Proposition \ref{=res}.
\hfill$\Box$

\subsection{Kato's residue homomorphisms for non-standard fields}\label{3.4}

For a two-dimensional local field $F$ of mixed characteristic which is not necessarily standard, choose a standard subfield $L$ inside $F$ such that
$F/L$ is finite and $k_F=k_L$. The symbol $\{,\}_F$ is defined to be the composition
\begin{equation}
\xymatrix{
F^\times\times F^\times\ar[r] & K_2(F) \ar[r]^{N_{F/L}} & K_2(L) \ar[r]^{\partial_L} & K_1(k_L),
}
\end{equation}
and we denote by $\partial_F$ the induced map $K_2(F)\to K_1(k_F)$. By \cite{K} Proposition 3, $\partial_F$ is independent of the choice of $L$.
If we let $e:=e(F/L)=e(F/k_L)$, then by \cite{K} Proposition 2 one has
\begin{equation}\label{ramify}
\partial_F(U^{ie}K_2(F))\subset \partial_L U^iK_2(L) \subset U^iK_1(k_F)=U^i_{k_F}.
\end{equation}
Hence we see that $\partial_F$ factors through $\widehat{K}_2(F)$. We also remark that $\partial_F$ induces a map $K_2({\mathcal O}_F)\to K_1({\mathcal O}_{k_F})$, similar with the case of standard fields.

\begin{prop}\label{funmix}
Let $F'/F$ be a finite extension of two-dimensional local fields of mixed characteristic. Then the following diagram commutes:
\begin{equation*}
\xymatrix{
K_2(F')\ar[r]^{\partial_{F'}} \ar[d]_{N_{F'/F}} & K_1(k_{F'})\ar[d]^{N_{k_{F'}/k_F}} \\
K_2(F) \ar[r]^{\partial_F} & K_1(k_F).
}
\end{equation*}
\end{prop}

\noindent\textit{Proof.}
Let $L$ be a standard subfield of $F$ used to define $\partial_F$. Then $L'=k_{F'}L$ is also standard and can be used to define
$\partial_{F'}$. By functoriality of norm homomorphisms for Milnor $K$-groups and Proposition \ref{funstd}, we have
\[
\partial_F N_{F'/F}=\partial_L N_{F'/L}=\partial_L N_{L'/L} N_{F'/L'}=N_{k_{L'}/k_L}\partial_{L'} N_{F'/L}=N_{k_{F'}/k_F}\partial_{F'}.
\]
\hfill$\Box$

\section{A reciprocity for two-dimensional normal local rings}\label{sec4}

Our settings in this section and the next are almost identical to those in \cite{Mo} and \cite{Mo2}, and in the proofs we shall apply many similar arguments
 as in \cite{Mo2} to the $K$-groups instead
of differential forms.

In this section, all rings are assumed to be Noetherian.

\subsection{Reciprocity for complete rings}\label{4.1}

Assume that $A$ is a two-dimensional normal complete local ring of characteristic zero with finite residue field of characteristic $p$. Let $F$
be the fraction field and $\mathfrak{m}$ be the maximal ideal of $A$. For each height 1 prime $y$ of $A$, the localization $A_y$ is a discrete valuation ring, and denote by $F_y$ the fraction field of $\widehat{A_y}$, which is a two-dimensional local field of characteristic zero.

Fix a finite extension ${\mathcal O}_k$ of $\mathbb{Z}_p$ inside $A$, where ${\mathcal O}_k$ is the ring of integers of a finite extension $k$ of $\mathbb{Q}_p$ inside $F$. For each height 1 prime $y$ of  $A$, the constant/coefficient field $k_y:=k_{F_y}$ of $F_y$ is a finite extension of $k$. Let $\{,\}_{F_y}: F_y^\times\times F_y^\times\to k_y^\times$
 be the symbol defined in section \ref{sec3}, and let $\{,\}_{F_y/k}$ be the composition $N_{k_y/k}\{,\}_{F_y}: F_y^\times\times F_y^\times \to k^\times$.

\begin{theorem}\label{comp}
Let $f,g\in F^\times$. Then $\{f,g\}_{F_y}=1$ for almost all height 1 primes $y$ of $A$, and in $k^\times$ one has
\[
\prod_{\mathrm{ht}(y)=1} \{f,g\}_{F_y/k}=1.
\]
\end{theorem}

\noindent\textit{Proof.}
By \cite{Mo} Lemma 3.7, there is a subring $B$ of $A$ which contains ${\mathcal O}_k$ and is ${\mathcal O}_k$-isomorphic to
${\mathcal O}_k[[T]]$, such that $A$ is a finite $B$-module. Then the theorem holds for $B$, by Proposition \ref{=res}
and Kato's reciprocity law (\cite{K} Proposition 4).
Let $M$ be the fraction field of $B$. Let us first prove that, for a fixed height 1 prime $y$ of $B$,
\begin{equation}\label{localglobal}
\partial_{M_y/k}j_{M_y/M}N_{F/M}\{f,g\}=\prod_{Y|y}\{f,g\}_{F_Y/k},
\end{equation}
where $\{f,g\}$ is the image of $(f,g)$ in $K_2(F)$, and $Y$ runs over the finitely many height 1 primes of $A$ lying over $y$.
Applying Proposition
\ref{lgnorm} for $B_y$ and $A_y:=(B\backslash y)^{-1}A$ one has the following local-global norm formula
\[
j_{M_{y}/M}\circ N_{F/M}=\prod_{Y|y}N_{F_{Y}/M_{y}}\circ j_{F_{Y}/F}: K_2(F)\to K_2(M_y).
\]
Composing $\partial_{M_y/k}$ on both sides of the last equation, we see that (\ref{localglobal}) follows from Propositions \ref{funequ}
and \ref{funmix}. For almost all height 1 primes $Y$ of $A$, we have $p\not\in Y$ and $f, g\in A_Y^\times$. In particular $F_Y$ is of zero residue characteristic, and in this case $\{f,g\}_{F_Y}=1$ by (\ref{classtame}). This proves $\{f,g\}_{F_Y}=1$ for almost all height one primes $Y$ of $A$. The theorem now follows from
(\ref{localglobal}) and the reciprocity for $B$.
\hfill$\Box$

\subsection{Reciprocity for incomplete rings}

Similar to \cite{Mo} section 3.3, we remark that it is possible to remove the restriction that $A$ be complete. Let ${\mathcal O}_k$ be a discrete valuation ring of characteristic zero with finite residue field. Assume that $A$ is a two-dimensional normal local ring with finite residue field, and that $A$ is the localization of a finitely-generated ${\mathcal O}_k$-algebra.

Let $\widehat{A}$ and $\widehat{\mathcal O}_k$ be the completions, which then satisfy the conditions in section \ref{4.1}. Let $F$, $\widehat{F}$ and $\hat{k}$ be the fraction fields of $A$, $\widehat{A}$ and $\widehat{\mathcal O}_k$ respectively. For a height 1 prime $y$ of $A$, we define
\begin{equation}
\{,\}_{F_y/\hat{k}}=\prod_{Y|y}\{,\}_{\widehat{F}_Y/\hat{k}}: F^\times\times F^\times \hookrightarrow \widehat{F}^\times \times \widehat{F}^\times\to \hat{k}^\times,
\end{equation}
where $Y$ runs over the finitely many height 1 primes of $\widehat{A}$ lying over $y$.

Following \cite{Mo} section 3.3, a prime $Y$ of $\widehat{A}$ is called \textit{transcendental} if $Y\cap A=0$. If the height 1 prime $Y$ of $ \widehat{A}$ is not transcendental,
then it is a prime minimal over $y\widehat{A}$ where $y=Y\cap A$. On the other hand, if $Y$ is transcendental, then $p\not\in Y$ and $F\subset \widehat{A}_Y$, which implies that $\{f,g\}_{\widehat{F}_Y/\hat{k}}=1$ for any $f,g\in F^\times\subset \widehat{A}_Y^\times$. Hence we deduce

\begin{theorem}\label{incomp}
Let $f,g\in F^\times$. Then $\{f,g\}_{F_y}=1$ for almost all height 1 primes $y$ of $A$, and in $\hat{k}^\times$ one has
\[
\prod_{\mathrm{ht}(y)=1} \{f,g\}_{F_y/\hat{k}}=1.
\]
\end{theorem}

\subsection{Geometrization: reciprocity around a point}\label{4.3}

Let ${\mathcal O}_K$ be a Dedekind domain of characteristic zero with finite residue fields, and let $K$ be its fraction field. Let $X$ be a two-dimensional normal scheme, flat and of finite type over $S=\textrm{Spec }{\mathcal O}_K$, with function field denoted by $K(X)$. Let $x\in X$ be a closed point lying over a closed point $s\in S$, and let $y\subset X$ be an irreducible curve (one-dimensional closed subscheme) passing through $x$. Then $A=\widehat{\mathcal O}_{X,x}$ satisfies the conditions in section \ref{4.1}, and contains the complete discrete valuation ring $\widehat{\mathcal O}_{K,s}$. Let $K(X)_x$ and $K_s$ be the fraction fields of $\widehat{\mathcal O}_{X,x}$ and ${\mathcal O}_{K_s}:=\widehat{\mathcal O}_{K,s}$ respectively.

Assume that the curve $y$ has formal branches $Y_1,\ldots, Y_n$
 at $x\in y$, i.e.
 \[
 y|_{\textrm{Spec }\widehat{\mathcal O}_{X,x}}=\bigcup_{1\leq i\leq n}Y_i,
 \]
where $Y_i$ is irreducible in $\textrm{Spec }\widehat{\mathcal O}_{X,x}$ for $1\leq i\leq n$. By abuse of notations, let the height 1 primes $y$ of ${\mathcal O}_{X,x}$ and $Y_i$ of  $\widehat{\mathcal O}_{X,x}$ be the local equations of $y$ and $Y_i$ at $x$ respectively. Then we have
\[
\{Y_i, i=1,\ldots, n\}=\{\textrm{height 1 primes }Y\textrm{ of }{\mathcal O}_{X,x}: Y|y\}.
\]
Let $\widehat{\mathcal O}_{x,Y}$ be the completion of the localization $(\widehat{\mathcal O}_{X,x})_Y$ of $\widehat{\mathcal O}_{X,x}$ with respect to $Y$, and let $K_{x,Y}$ be its fraction field, which is a two-dimensional local field. Then we define
\begin{equation}\label{twod}
K_{x,y}:=\prod_{Y|y}K_{x,Y},\quad \widehat{\mathcal O}_{x,y}:=\prod_{Y|y}\widehat{\mathcal O}_{x, Y}.
\end{equation}
 Notice that there are embeddings $K(X)\hookrightarrow K(X)_x\hookrightarrow K_{x,y}$. Now we let
\begin{equation}\label{twodtame}
\{,\}_{x,y}\stackrel{\textrm{def}}{=}\prod_{Y|y}\{,\}_{K_{x,Y}/K_s}: K_{x,y}^\times \times K_{x,y}^\times\to K_s^\times,
\end{equation}
where $\{,\}_{K_{x,Y}/K_s}: K_{x,Y}^\times \times K_{x,Y}^\times\to K_s^\times$ is as defined in section \ref{sec3}. Then we have the following reciprocity around a point, which is
a geometric translation of Theorem \ref{incomp}.

\begin{theorem} \label{recipt}
Fix a closed point $x$ and let $f,g\in K(X)^\times$. Then
$\{f,g\}_{x,y}=1$ for almost all $y$ passing through $x$, and in $K_s^\times$ one has
\[
\prod_{y\subset X, y\ni x} \{f,g\}_{x,y}=1.
\]
\end{theorem}

\noindent\textit{Proof.}
By Theorem \ref{comp} one has
\[
\prod_{\textrm{ height 1 primes }Y\textrm{ of }\widehat{\mathcal O}_{X,x}}\{f,g\}_{K_{x,Y}/K_s}=1.
\]
If a height 1 prime $Y$ of $\widehat{\mathcal O}_{X,x}$ is not a formal branch at $x$ of some irreducible curve $y\subset X$, then $Y$ is transcendental and
$\{f,g\}_{K_{x,Y}/K_s}=1$. Hence the theorem follows.
\hfill$\Box$

\section{Reciprocities for arithmetic surfaces}\label{sec5}

Let ${\mathcal O}_K$ and $K$ be as in section \ref{4.3}, and let $X$ be an ${\mathcal O}_K$-curve, i.e. a normal scheme, proper and flat over $S=\textrm{Spec }{\mathcal O}_K$, whose generic fibre $X_K$ is a smooth and geometrically
connected curve. We shall keep the notations in section \ref{4.3}. Moreover, for an irreducible curve $y\subset X$, let $K(X)_y$ be the fraction field of $\widehat{\mathcal O}_{X,y}$, the completion of the discrete valuation ring ${\mathcal O}_{X,y}$. Then there are embeddings $K(X)\hookrightarrow K(X)_y\hookrightarrow K(X)_{x,y}$, where $x\in y$ is a closed point.

\subsection{Reciprocity along vertical curves}\label{5.1}

We shall establish the following reciprocity law for vertical curves on an arithmetic surface, which is similar to \cite{Mo2} Theorem 3.1.

\begin{theorem}\label{recivertical}
Let $y\subset X$ be an irreducible component of a special fibre $X_s$, where $s\in S$ is a closed point. Let $f,g\in K(X)_y^\times$. Then
$\prod\limits_{x\in y}\{f,g\}_{x,y}$ converges to $1$ in $K_s^\times$, where the product is taken over all closed points $x$ of $y$.
\end{theorem}

As preparations, let us first prove the convergence of the product in the last theorem (cf. \cite{Mo} Lemma 3.3).

\begin{lemma}\label{convergence}
With the conditions in Theorem \ref{recivertical}, $\prod\limits_{x\in y}\{f,g\}_{x,y}$ converges in $K_s^\times$. Moreover, the pairing
\[
K(X)_y^\times\times K(X)_y^\times \to K_s^\times, \quad (f,g)\mapsto \prod_{x\in y}\{f,g\}_{x,y}
\]
is continuous.
\end{lemma}

\noindent\textit{Proof.}
Let $\pi$ be a parameter of $K_s$ and let $e_y=\nu_y(\pi)$, where $\nu_y$ is the normalized valuation on $K(X)_y$. Since $K(X)$ is dense in $K(X)_y$, for any $n\geq 1$ we may write $f=f_n u_n$, $g=g_n v_n$ with $f_n, g_n\in K(X)^\times$, $u_n, v_n\in U^{ne_y}_{K(X)_y}$. Almost all $x\in y$ satisfy the following conditions:

(i) $X_s$ has no other irreducible components passing through $x$,

(ii) $x\not\in y'$ for any horizontal $y'$ appearing in $\textrm{div}(f_n)$ or $\textrm{div}(g_n)$.
\\
For such $x$, by Theorem \ref{recipt} we have
\[
\{f_n, g_n\}_{x,y}=\prod_{y'\ni x, y'\textrm{ horizontal}} \{f_n,g_n\}^{-1}_{x,y'}=1,
\]
where the last equality follows because $f_n, g_n\in \mathcal O_{X,y'}^\times$ and $K_{x,y'}$ is of equal characteristic for each horizonal $y'\ni x$.
Then from (\ref{ramify}) we deduce that
\[
\{f,g\}_{x,y}=\{f_n,v_n\}_{x,y}\{u_n,g\}_{x,y}\in U^n_{K_s}.
\]
Hence we have proved that for any $n\geq1$, $\{f,g\}_{x,y}\in U^n_{K_s}$ for almost all $x\in y$. This implies $\prod\limits_{x\in y}\{f,g\}_{x,y}$ converges in $K_s^\times$. Guaranteeing the convergence, we have $\prod\limits_{x\in y}\{f,g\}_{x,y}\in U^n_{K_s}$ whenever $f$ or $g$ lies in $U^{ne_y}_{K(X)_y}$, again using (\ref{ramify}). Thus we obtain the continuity of the pairing.
\hfill$\Box$

\

We also need the following simple but useful lemma, which is similar to \cite{Mo2} Lemma 5.1.

\begin{lemma}\label{generic}
Let $C$ be a smooth and geometrically connected curve over a field $K$ of characteristic zero, $L$ an extension of $K$, and $z$ a closed point of
$C$.

(i) Let $z'$ be a closed point of $C_L$ lying over $z$. Then the following diagram commutes:
\begin{equation*}
\xymatrix{
K(C)^\times_z\times K(C)_z^\times \ar[r]^{\quad\qquad\{,\}_z} \ar[d] & k(z)^\times\ar[d]\\
K(C_L)_{z'}^\times\times K(C_L)_{z'}^\times \ar[r]^{\quad\qquad\{,\}_{z'}} & k(z')^\times,
}
\end{equation*}
where $\{,\}_z$ is the one-dimensional tame symbol associated to the closed point $z$ on $C$, and $k(z)$ is the residue field of $z$; $\{,\}_{z'}$ and
$k(z')$ are defined similarly.

(ii) Let $z'$ run over all closed points of $C_L$ lying over $z$. Then the following diagram commutes:
\begin{equation*}
\xymatrix{
K(C)^\times_z\times K(C)_z^\times \ar[rr]^{\qquad N_{k(z)/K}\{,\}_z} \ar[d] && K^\times\ar[d]\\
\prod\limits_{z'|z}(K(C_L)_{z'}^\times\times K(C_L)_{z'}^\times) \ar[rr]^{\qquad\qquad\prod\limits_{z'|z}N_{k(z')/L}\{,\}_{z'}} && L^\times.
}
\end{equation*}
\end{lemma}

\noindent\textit{Proof.}
(i) follows from the compatible isomorphisms $K(C)_z\cong k(z)((t))$, $K(C_L)_{z'}\cong k(z')((t))$, where $t\in K(C)$ is a local parameter at $z$.
Since $k(z)\otimes_K L\cong \prod\limits_{z'|z}k(z')$, we have $N_{k(z)/K}=\prod\limits_{z'|z}N_{k(z')/L}$. (ii) follows from this fact together with (i).
\hfill$\Box$

\

Let $X'=X\times_{{\mathcal O}_K} {\mathcal O}_{K_s}$, and $p: X'\to X$ be the natural morphism. Then $p$ induces
an isomorphism of special fibres
$X'_s\cong X_s$, and for any $x'\in X_s'$, $p$ induces an isomorphism of complete local rings $\widehat{\mathcal O}_{X,p(x')}\cong\widehat{\mathcal O}_{X',x'}$. Note that $X'$ is an ${\mathcal O}_{K_s}$-curve with generic fibre $X_{K_s}$. The following lemma is similar to
\cite{Mo} Lemma 3.5.

\begin{lemma}\label{completion}
Let $y\subset X$ be an irreducible curve and $x\in y$ be a closed point lying over $s$. Then the following diagram commutes:
\begin{equation*}
\xymatrix{
\prod\limits_{y'|y, y'\ni x'}(K(X')_{y'}^\times\times K(X')_{y'}^\times) \ar[drr]^{\prod_{y'|y}\{,\}_{x',y'}}  && \\
K(X)_y^\times\times K(X)_y^\times \ar[u] \ar[rr]^{\{,\}_{x,y}} && K_s^\times,
}
\end{equation*}
where $x'$ is the unique closed point of $X'_s$ lying over $x$, and $y'$ runs over irreducible curves of $X'$ which lie over $y$ and contain $x'$.
\end{lemma}

\noindent\textit{Proof.}
This follows from definitions (\ref{twod}) and (\ref{twodtame}). Indeed, by abuse of notation let the height 1 prime $y$ of ${\mathcal O}_{X,x}$ be the local equation of $y$. Then
\[
\{,\}_{x,y}=\prod_{\textrm{height 1 primes }Y\textrm{ of }\widehat{\mathcal O}_{X,x}, Y|y}\{,\}_{K_{x,Y}/K_s}.
\]
Similarly, one has
\begin{eqnarray*}
\prod_{y'|y, y'\ni x'}\{,\}_{x',y'}&=&\prod_{\textrm{height 1 primes }y'\textrm{ of }{\mathcal O}_{X',x'}, y'|y}\left[\prod_{\textrm{height 1 primes }Y'\textrm{ of }\widehat{\mathcal O}_{X',x'}, Y'|y'}\{,\}_{K_{x',Y'}/K_s}\right]\\
&=&\prod_{\textrm{height 1 primes }Y'\textrm{ of }\widehat{\mathcal O}_{X',x'}, Y'|y}\{,\}_{K_{x',Y'}/K_s}.
\end{eqnarray*}
The equality of maps $\{,\}_{x,y}=\prod\limits_{y'|y, y'\ni x'}\{,\}_{x',y'}$ follows from the isomorphism $\widehat{\mathcal O}_{X,x}\cong \widehat{\mathcal O}_{X',x'}$ remarked above.
\hfill$\Box$

\

Let $z$ be the generic point of $y$. Since ${\mathcal O}_{X, y}={\mathcal O}_{X, z}\cong {\mathcal O}_{X_K, z}$ (see e.g. \cite{L} Chap 8. Exercise 3.7) and $X_K$ has the same function field as
 $X$, we see that $K(X)_y=K(X_K)_z$.

\begin{cor}\label{one-two}
Let $y$ be an irreducible horizontal curve on $X$, whose generic point $z$ is a closed point of $X_K$. Then for $f, g\in K(X)_y^\times$ one has
\[
\prod_{x\in y\cap X_s}\{f,g\}_{x,y}=N_{k(z)/K}\{f,g\}_z
\]
where $\{,\}_z$ is the usual tame symbol for the curve $X_K$ at the closed point $z$.
\end{cor}

\noindent\textit{Proof.}
Applying Lemma \ref{generic} for $C=X_K$ and $L=K_s$ we obtain
\[
N_{k(z)/K}\{f,g\}_z=\prod_{z'|z}N_{k(z')/K_s}\{f,g\}_{z'},
\]
where $z'$ runs over all closed points of $X_{K_s}$ lying over $z$. Each closed point $z'\in X_{K_s}$ has a unique reduction $x'$ on $X'_s$, i.e. $\overline{z'}$ meets the special fibre $X'_s$ at a unique point $x'$. This implies that the irreducible curve $y':=\overline{z'}$ on $X'$ sits over $y=\bar{z}$ and contains $x'$. Therefore by Lemma \ref{completion} we obtain $\{f,g\}_{x,y}=\{f,g\}_{x',y'}$, where $x=p(x')\in y\cap X_s$. To finish the proof it remains to show that
$N_{k(z')/K_s}\{f,g\}_{z'}=\{f,g\}_{x',y'}$.

It can be shown that (see e.g. \cite{Mo2} Lemma 3.8) that $\widehat{\mathcal O}_{X_{K_s},z'}=\widehat{\mathcal O}_{x',y'}$, where $\widehat{\mathcal O}_{x',y'}$ is the completion of
$\widehat{\mathcal O}_{X',x'}$ with respect to the discrete valuation given by the height 1 prime $y'$ of  $\widehat{\mathcal O}_{X',x'}$. Then the one-dimensional tame symbol for the curve $X_{K_s}$ at $z'$ gives the tame symbol for the two-dimensional local field $K_{x',y'}=\textrm{Frac }\widehat{\mathcal O}_{X_{K_s},z'}$, which is of equal characteristic.
\hfill$\Box$

\

Using arguments similar to the proof of \cite{Mo2} Theorem 3.1, now we can prove Theorem \ref{recivertical}. Let $y_1:=y, y_2, \ldots, y_l$ be the irreducible components of $X_s$. By the continuity proved in Lemma \ref{convergence} we only need to show that $\prod\limits_{x\in y}\{f,g\}_{x,y}=1$ for any $f, g\in K(X)^\times$. From the Weil reciprocity law for the curve
$X_{K}$ and Corollary \ref{one-two}, it follows that
\[
\prod_{y\textrm{ horizontal}}\prod_{x\in y\cap X_s}\{f,g\}_{x,y}=1.
\]
Since $\{f,g\}_{z'}=1$ for almost all closed points $z'\in X_{K_s}$, from the proof of previous corollary we see that
only finitely many terms in the last product differ from 1. Hence we may rewrite it as
\[
\prod_{x\in X_s}\prod_{y\ni x, y\textrm{ horizontal}}\{f,g\}_{x,y}=1.
\]
Using Theorem \ref{recipt} we deduce that
\[
\prod_{x\in X_s}\prod_{y\ni x, y\textrm{ vertical}}\{f,g\}_{x,y}=1,
\]
which by Lemma \ref{convergence} can be rearranged as
\[
\prod^l_{i=1}\prod_{x\in y_i}\{f,g\}_{x,y_i}=1.
\]
Let $\nu_i$ be the discrete valuation on $K(X)$ corresponding to $y_i$, $i=1,\ldots, l$. For any $n\geq 1$ choose $f_n\in K(X)^\times$ such that
$\nu_1(f_n-1)\geq n$, $\nu_i(f_n-f)\geq n$, $i=2,\ldots, l$. Replacing $f$ by $f/f_n$ and letting $n\to \infty$, we see that
\[
\prod_{x\in y}\{f,g\}_{x,y}=\lim_{n\to\infty}\prod^l_{i=1}\prod_{x\in y_i}\{f/f_n,g\}_{x,y_i}=1.
\]
Here we have used Lemma \ref{convergence}. This finishes the proof of Theorem \ref{recivertical}.

\subsection{Reciprocity along horizontal curves}\label{5.2}

We shall follow the treatment in \cite{Mo2} section 5. Assume that ${\mathcal O}_K$ is the ring of integers of a number field $K$. Let $\overline{S}=S_f\cup S_\infty$ denote the set of places of $K$, where $S_f$ denotes the finite places, i.e. closed points of $S=\textrm{Spec } {\mathcal O}_K$, and $S_\infty$ the infinite ones. Write $X_v$ for $X_{K_v}$, where $v\in\overline{S}$. Note that for $v\in S_\infty$, $X_v$ is a smooth projective curve over $K_v=\mathbb{R}$ or $\mathbb{C}$.

Let $y$ be an irreducible horizontal curve on $X$. Then its generic point $z$ is a closed point of $X_K$. For any place $v$, let $y\cap X_v$ denote the set of closed points of $X_v$ lying over $z$, i.e. the preimage of $z$ under the projection $X_{v}\to X_K$. Let $k(z)$ be the residue field of $z$, which is a finite extension of $K$. Then $y\cap X_v$ is finite and corresponds to the places of $k(z)$ extending $v$. We now compactify $y$ by adding all the finite sets $y\cap X_v$ for $v\in S_\infty$, and by abuse of notation still denote by $y$ the compactification.

For $x\in y\cap X_v$, define a symbol $\{,\}_{x,y}: K(X)_y^\times\times K(X)_y^\times\to K_v^\times$ by
\begin{equation*}
\xymatrix{
K(X)_y^\times\times K(X)_y^\times\ar[r] & K(X_v)_x^\times\times K(X_v)_x^\times \ar[r]^{\quad\qquad\{,\}_x} & k(x)^\times\ar[r]^{N_{k(x)/K_v}} & K_v^\times,
}
\end{equation*}
where $K(X_v)$ is the function field of $X_v$, $\{,\}_x$ is the one-dimensional tame symbol associated to the closed point $x$ on $X_v$, and $k(x)$ is the residue field of $x$.

Applying Lemma \ref{generic} for $C=X_K$ and $L=K_v$ we obtain

\begin{cor}\label{glob}
With previous notations, for $f, g\in K(X)_y^\times$ one has
\[
\prod_{x\in y\cap X_v}\{f,g\}_{x,y}=N_{k(z)/K}\{f,g\}_z.
\]
\end{cor}

Let $\mathbb{I}_K$ be the idele group of $K$, and $\chi=\bigotimes\limits_{v\in\overline{S}}\chi_v: \mathbb{I}_K \to S^1$ be an idele class character, i.e. a continuous character of $\mathbb{I}_K$ which is trivial on $K^\times$. For a closed point $x\in y$ lying over $v\in\overline{S}$, define the symbol $\chi_{x,y}$ to be the composition
\begin{equation*}
\xymatrix{
\chi_{x,y}: K(X)_y^\times\times K(X)_y^\times \ar[r]^{\qquad\qquad\{,\}_{x,y}} & K_v^\times \ar[r]^{\chi_v} & S^1.
}
\end{equation*}
The following reciprocity along horizontal curves is analogous to \cite{Mo2} Theorem 5.4.

\begin{theorem}\label{recihori}
Assume that $\mathcal{O}_K$ is the ring of integers of a number field $K$. Let $y$ be an irreducible horizontal curve on $X$ and $f, g\in K(X)_y^\times$. Then $\chi_{x,y}(f,g)=1$ for almost all closed points $x\in y$, and in $S^1$ one has
\[
\prod_{x\in y}\chi_{x,y}(f,g)=1.
\]
\end{theorem}

\noindent\textit{Proof.}
Let $\pi$ be the morphism $X\to S$. Using that $K(X)$ is dense in $K(X)_y$, we may write $f=\tilde{f}u_1$, $g=\tilde{g}u_2$, where
$\tilde{f}, \tilde{g}\in K(X)^\times$, $u_1, u_2\in U^1_{K(X)_y}$. Then almost all $x\in y$ satisfy the condition that
$x\not\in y'$ for any irreducible curve $y'\neq y$ which appears in $\textrm{div}(\tilde{f})$ or $\textrm{div}(\tilde{g})$.
For such $x$, we have $\tilde{f},\tilde{g}\in\widehat{\mathcal O}_{X,y'}^\times$ for any $y'\ni x$, $y'\neq y$. Then by Theorem \ref{recipt} we have
\[
\{\tilde{f},\tilde{g}\}_{x,y}=\prod_{y'\ni x, y'\neq y}\{\tilde{f},\tilde{g}\}^{-1}_{x,y'}=\prod_{y'\ni x, y'\textrm{ vertical}}\{\tilde{f},\tilde{g}\}^{-1}_{x,y'}\in {\mathcal O}_{K_{\pi(x)}}^\times.
\]
By properties of tame symbols for two-dimensional local fields of equal characteristic, we also have
\[
\{\tilde{f},u_2\}_{x,y}=\{u_1,\tilde{g}\}_{x,y}=\{u_1,u_2\}_{x,y}=1
\]
for any $x\in y$. Therefore we have proved that $\{f,g\}_{x,y}\in \mathcal O^\times_{K_{\pi(x)}}$ for almost all $x\in y$. Since
the character $\chi_v$ is unramified for almost all $v\in S_f$, we obtain the first statement of the theorem. The idele
\[
\left(\prod_{x\in y\cap X_v}\{f,g\}_{x,y}\right)_{v\in\overline{S}}
\]
is global by Corollary \ref{one-two} and \ref{glob}, hence the product formula follows.
\hfill$\Box$

\section*{Acknowledgement}

This work was started during the author's visit at the Hong Kong University of Science and Technology, and partially supported by NSFC11201384. The author would like to heartily thank Prof. Y. Zhu
for helpful discussions. The author is also grateful to the referee for many valuable comments on a preliminary version of this paper.

\end{document}